\documentclass[11pt,a4paper]{article}
\usepackage{latexsym}
\usepackage{amsfonts}
\usepackage{amsmath}
\usepackage{amssymb}
\usepackage{graphicx}
\numberwithin{equation}{section}
\newtheorem{theorem}{Theorem}[section]

\newtheorem{notation}[theorem]{Notation}
\newtheorem{corollary}[theorem]{Corollary}
\newtheorem{definition}[theorem]{Definition}

\newtheorem{lemma}[theorem]{Lemma}
\newtheorem{proposition}[theorem]{Proposition}
\newtheorem{remark}[theorem]{Remark}
\newtheorem{conj}[theorem]{Conjecture}

\begin{document}
\title{Beltrami operators,\\ non--symmetric elliptic equations\\ and quantitative Jacobian bounds}
\author{Giovanni Alessandrini\\
\small Dipartimento di Matematica e Informatica,
Universit\`a  di Trieste\\
\small Via Valerio 12/b, 34100 Trieste, Italia, e-mail:
\texttt{alessang@units.it}\\
\vspace*{0.05cm}\\
Vincenzo Nesi \\
\small
Dipartimento di Matematica,
La Sapienza, Universit\`a di Roma,\\
\small P. le A. Moro 2, 00185 Roma, Italia, e-mail:
\texttt{nesi@mat.uniroma1.it}}
\date{\today}

\maketitle

\begin{abstract}
\noindent In recent studies on the G-convergence of Beltrami operators, a number
of issues arouse concerning injectivity properties of families of quasiconformal
mappings.
Bojarski, D'Onofrio, Iwaniec and Sbordone formulated a conjecture based on the existence of a so-called primary pair.
Very recently, Bojarski proved the existence of one such pair. We provide a general, constructive, procedure for obtaining a new rich class of
such primary pairs.

This proof is obtained as  a slight
adaptation of previous work by the authors concerning the
nonvanishing of the Jacobian of pairs of solutions of elliptic
equations in divergence form in the plane. It is proven here that
the results previously obtained when the coefficient matrix is
symmetric also extend to the non-symmetric case. We also prove a
much stronger result giving a quantitative bound for the Jacobian
determinant of the so-called \emph{periodic} $\sigma$-harmonic
sense preserving homeomorphisms of $\mathbb C$ onto itself.
\end{abstract}

\noindent{\small \textit{2000 AMS Mathematics Classification Numbers}: 30C62, 35J55}

\noindent{\small \textit{Keywords}:  Beltrami operators, quasiconformal mappings}

\section{Introduction.}
In order to explain the results of this paper and their
motivations, it is necessary to introduce a number of topics, and
to illustrate their mutual relationships. These topics are
Beltrami operators and their  associated concept of $G$-convergence, non-symmetric elliptic
operators in divergence form and $H$-convergence,
$\sigma$-harmonic mappings.

\subsection{The $G$-convergence of Beltrami operators and the $K>3$ conjecture.}
Recently Iwaniec et al. \cite{B0} and Bojarski et al. \cite{B},
introduced a notion of $G$-convergence for Beltrami operators,
aimed at generalizing to this context the well-known theory of
$G$-convergence initiated by Spagnolo \cite{spa} and De Giorgi
\cite{despa}.  Let us recall their definitions and the main
conjecture in \cite{B}. Let $\Omega$ be a bounded, simply
connected open subset of $\mathbb R^2$, and, as usual, let us
identify points $x=(x_1,x_2)\in\mathbb R^2$ with points $z\in
\mathbb C$ through the relation $z=x_1+ix_2$ .   Let $\nu$ and
$\mu$ be two complex valued measurable functions defined on
$\Omega$ and satisfying, for some $K\geq 1$, the following
ellipticity condition
\begin{equation}\label{ellQC}
|\mu|+|\nu|\leq \frac{K-1}{K+1} \ .
\end{equation}
Consider the following first order non homogeneous Beltrami equation
\begin{equation}\label{1stordernh}
\begin{array}{ll}
f_{\overline{z}}-\mu f_z -\nu \overline{f_z}= g \ .
\end{array}
\end{equation}
Given a sequence of pairs of Beltrami coefficients $(\mu_j,
\nu_j)$ and an extra pair $(\mu,\nu)$ all satisfying \eqref{ellQC}, for a fixed  $K\geq 1$, one
denotes by $\mathcal B_j$,  $\mathcal B$ the differential operators defined  as
follows
\begin{equation}\label{boj}
\mathcal B_j := \frac{\partial}{ \partial \overline{z}}-\mu_j
\frac{\partial}{ \partial z} -\nu_j \overline{\frac{\partial}{
\partial z}} \ ,
\end{equation}
\begin{equation}\label{bo}
\mathcal B := \frac{\partial}{ \partial \overline{z}}-\mu \frac{\partial}{ \partial z} -\nu \overline{\frac{\partial}{ \partial z}} \ ,
\end{equation}
so that (\ref{1stordernh})  can be rewritten as
\begin{equation}
\mathcal B f= g \ .
\end{equation}
The authors in  \cite{B0} introduce the following definition, and prove Theorem~\ref{main-B} below.
\begin{definition}
The sequence of differential operators $ \mathcal B_j$ is said to
$G$-converge to $ \mathcal B$ if, for any sequence $f_j\in
W^{1,2}(\Omega;\mathbb C)$ which converges weakly to $f\in
W^{1,2}(\Omega;\mathbb C)$,    and such that  $\mathcal B_j f_j$ converges strongly in $L^2(\Omega;\mathbb C)$, one has
\begin{equation}
\lim_{j\to +\infty} \mathcal B_j f_j =\mathcal B f
\end{equation}
strongly in $L^2(\Omega;\mathbb C)$.
\end{definition}
\begin{theorem}[\cite{B0}]\label{main-B}
For any $K\in [1,3]$, the family of Beltrami operators defined by
(\ref{bo}) and satisfying (\ref{ellQC}) is $G$-compact.
\end{theorem}

In order to explain our new main results and to put the previous
one into context, let us begin by explaining the main point in the
proof of Theorem~\ref{main-B}.

As previously outlined one of the main results in \cite{B0}  is a
compactness result obtained
 under an assumption of  \emph{small ellipticity}, that is,  $K\leq 3$ in \eqref{ellQC}.

The key to this result relies on the following issue. Let $\Omega$
be a bounded, open and convex set. Let $(\mu,\nu)$ be a Beltrami
pair satisfying (\ref{ellQC}) and let $\Phi$ and $\Psi$ be the
solutions to
 \begin{equation}\label{pair1storder}
 \left\{
\begin{array}{lr}
\Phi_{\bar{z}}=\mu \Phi_z +\nu \overline{\Phi_z}\ , & \hbox{in $\Omega$}\ ,\\
\mathfrak{Re} \Phi = x_1 \ , & \hbox{on $\partial \Omega$}\ , \\
\Psi_{\bar{z}}=\mu \Psi_z +\nu \overline{\Psi_z}\ ,& \hbox{in $\Omega$}\ ,\\
\mathfrak{Re} \Psi = x_2\ , & \hbox{on $\partial \Omega$} \ ,
\end{array}
\right.
\end{equation}
where the boundary conditions are  understood in
the sense of  $W^{1,2}(\Omega)$ traces. The pair $(\Phi,\Psi)$ is called a primary pair.
In \cite{B} the authors formulate the following conjecture.

\begin{conj}\label{conj1}
Let $(\mu,\nu)$ be complex valued measurable  coefficients
satisfying (\ref{ellQC}). Then the pair of quasiconformal mappings
$\Phi$ and $\Psi$ defined by (\ref{pair1storder}) satisfies the
following pointwise inequality:
 \begin{equation}\label{det>0C}
\begin{array}{llll}
\mathfrak{Im} (\Phi_z \overline{\Psi_z})>0&\hbox{almost everywhere}&\hbox{in}&\Omega.
\end{array}
\end{equation}
\end{conj}
In \cite{B0,B} it is proven that, if Conjecture \ref{conj1} holds, then Theorem \ref{main} follows.

As a consequence of our results we prove that \eqref{det>0C} holds and therefore we obtain  the following result.
\begin{theorem} \label{main}
For any $K\in [1,+\infty)$, the family of Beltrami operators
defined by (\ref{bo}) and satisfying (\ref{ellQC}) is $G$-compact.
\end{theorem}

Very recently, Bojarski \cite{bo2} has proved a result which also implies Theorem \ref{main} but does not solve Conjecture \eqref{conj1}. More precisely he has proven that given $\Omega$ and a Beltrami pair $(\mu,\nu)$ satisfying \eqref{ellQC} there exists a   {\em primary pair} $(\Phi,\Psi)$ so that $\Phi$ and $\Psi$  are quasiconformal mappings of the complex plane onto itself satisfying the Beltrami equations with coefficients $\mu$ and $\nu$ and satisfy \eqref{det>0C}. Bojarski's primary pair is obtained by requiring the so-called {\em hydrodynamical normalization}, that is, by looking for a globally homeomorphic solution of $\mathbb C$ onto itself obtained as follows. First extend  $(\mu,\nu)$ to be zero in the complement of $\Omega$. Then look for a solution of the new Beltrami equation defined on $\mathbb C$. Such a solution will be holomorphic near infinity. Then normalize the behaviour at  infinity of such function. By the seminal work of Bojarski (see the references of \cite{bo2}), it is known that  one obtains a quasiconformal mapping of $\mathbb C$ onto itself. This beautiful construction however does not set the question of whether   the Dirichlet data in \eqref{pair1storder} will provide us with a primary pair. We prove that this is the case in Theorem \ref{conj2}. In fact we provide a large class of Dirichlet boundary data achieving the desired task. We use the combination of Theorem \ref{Kneser}
and Theorem \ref{t3.1}. See Corollary \ref{pr-sol}.

\subsection{Second order equations in divergence form, ellipticity and $H$-convergence.}
It is well known that Beltrami equations with complex dilatations
$\nu$ and $\mu$ give rise in a very natural way to second order
elliptic operators whose coefficient matrices $\sigma$ depend in
an explicit way upon $\nu$ and $\mu$ and conversely. A brief
review will be offered in the following subsection.
The authors in \cite{B0, B} use the notion of $G$-convergence for
Beltrami operators also to induce a concept of $G$-convergence for
second order non-symmetric operators in divergence form (see
Definition 2 in \cite{B}) and to treat the $G$-convergence of
second order non-divergence equations (see \cite{B0}). We shall
not enter such issues in this note, however we observe that it is also
instructive to recall the notion $H$-convergence introduced by
Murat and Tartar for possibly non-symmetric, elliptic operators
in divergence form. An easily accessible reference is \cite{mt}. The original work dates back to 1977 (see the quoted reference for more details).

\begin{definition}\label{def.H}
Consider a bounded, open, simply connected set $\Omega\subset
\mathbb R^2$. Given positive constants $\alpha$ and $\beta$, we
say that a measurable function $\sigma$, defined on $\Omega$ with
values into the space of $2\times2$ matrices,  belongs to the
class $\mathcal M (\alpha,\beta,\Omega)$  if one has
\begin{equation}\label{ellTM}
\begin{array}{ccrllll}
\sigma(z) \xi\cdot \xi&\geq& \alpha |\xi|^2&,&\hbox{for every $\xi \in \mathbb R^2$ and for a.e. $z\in\Omega$ ,}\\
\sigma^{-1}(z) \xi \cdot \xi &\geq & \beta^{-1}
|\xi|^2&,&\hbox{for every $\xi \in \mathbb R^2$ and for a.e.
$z\in\Omega$ \ .}
\end{array}
\end{equation}
\end{definition}
It is obvious that, for $\lambda=\alpha$ and for some $M > 0$,
such bounds are equivalent to the usual ellipticity bounds for
second order elliptic operators, see for instance \cite[Chapter
8]{gt}
\begin{equation}\label{ellGT}
\begin{array}{ccrllll}
\sigma(z) \xi\cdot \xi&\geq& \lambda |\xi|^2,&&\hbox{for
every $\xi \in \mathbb R^2$ and for a.e. $z\in\Omega$ ,}\\
\sum\limits_{i,j=1}^2|\sigma_{ij}(z)|^2& \leq & M \ ,& & \hbox{for
a.e. $z\in\Omega$ .}
\end{array}
\end{equation}
Yet another notion, originally used for the $H$-convergence is the following.
\begin{definition}
A matrix $\sigma$ with measurable entries belongs to $M(\lambda,\Lambda,\Omega)$ if
\begin{equation}\label{ellTMold}
\begin{array}{cclllll}
\sigma(z) \xi\cdot \xi&\geq& \lambda |\xi|^2\ ,&&\hbox{for
every $\xi \in \mathbb R^2$ and for a.e. $z\in\Omega$ ,}\\
|\sigma(z) \xi |& \leq & \Lambda|\xi|\ ,&&\hbox{for every $\xi \in
\mathbb R^2$ and for a.e. $z\in\Omega$ .}
\end{array}
\end{equation}
\end{definition}
However, different ways of bounding sets of matrices $\sigma$ may
or may not give rise to compact classes with respect to
convergences of weak type. To explain this let us recall the notion of $H$-convergence \cite{mt}.

\begin{definition}
We say that a sequence of elliptic matrices $\sigma_{j} \in
\mathcal M (\alpha,\beta,\Omega)$ $H$-converges to $\sigma_0 \in
\mathcal M (\alpha,\beta,\Omega)$ if for any  $f\in
H^{-1}(\Omega)$ the weak solution $u_j$ to
\begin{equation}\label{2ndorder}
\begin{array}{llll}
-{\rm div}(\sigma_{j} \nabla u_{j} )=f \ ,& \hbox{in $\Omega$}\ ,&
u_{j}\in W^{1,2}_0(\Omega)\ ,
\end{array}
\end{equation}
satisfies the following properties
\begin{equation}\label{weakconv}
\left\{
\begin{array}{llll}
u_{j} \rightharpoonup u_0\ , &\hbox{weakly in $W^{1,2}(\Omega)$}\ ,\\
\sigma_{j} \nabla u_{j}\rightharpoonup \sigma_0 \nabla u_0 \ ,
&\hbox{weakly in $L^{2}(\Omega)$ ,}
\end{array}
\right .
\end{equation}
where $u_0$ denotes the weak solution to
\begin{equation}\label{2ndorder0}
\begin{array}{llll}
-{\rm div}(\sigma \nabla u_{0} )=f \ , & \hbox{in $\Omega$ ,}&
u_{0}\in W^{1,2}_0(\Omega) \ .
\end{array}\end{equation}
\end{definition}
One of the main results in this theory is compactness. Given any
sequence $\{\sigma_{j}\}\subset \mathcal M (\alpha,\beta,\Omega)$
there exists a subsequence which $H$-converges to some element of
${\mathcal M} (\alpha,\beta,\Omega)$. It is worth noting here
 that the compactness does indeed depend on the specific character of the ellipticity bounds given by Murat and Tartar.
 For instance, it is known that the set of matrices in $M(\lambda,\Lambda,\Omega)$, that is the set  constrained by \eqref{ellTMold},  is \emph{not} compact for $H$-convergence. Murat and Tartar proved that a sequence of matrices in
 $M(\lambda,\Lambda,\Omega)$ admits  (up to subsequence) an $H$-limit in the class $M\left(\lambda,\frac{\Lambda^2}{\lambda},\Omega\right)$.
 An explicit example given by Marcellini in \cite{mar} shows that there exist a sequence $\{\sigma_j\}\subset  M(\lambda,\Lambda,\Omega)$ such that
  its $H$-limit $\sigma_0$ is constant (with respect to position) and
  satisfies
 \begin{equation*}
\inf_{|\xi|=1}\sigma_0\xi\cdot\xi=\lambda \ ,
\sup_{|\xi|=1}|\sigma_0\xi| =(\Lambda^2/\lambda) \ .
\end{equation*}
Let us also recall that the approach of Murat and Tartar has been
later extended to larger classes of operators (under the name of
$G$-convergence) by Dal Maso, Chiad\`o-Piat and Defranceschi
\cite{dmcpdf}.

\subsection{Beltrami equations, second order equations in divergence form and ellipticity.}
Let us recall now the basic algebraic relationship between second order elliptic equations in divergence form and linear first order systems.
Given $\sigma \in {\mathcal M}(\alpha, \beta,\Omega)$,   let $u\in W^{1,2}_{\rm loc}(\Omega)$ be a weak solution to
\begin{equation}\label{e2.2}
\begin{array}{ll}
{\rm div}(\sigma \nabla u )=0 & \hbox{in $\Omega$}\ .
\end{array}
\end{equation}
Then there exists $\tilde u \in W^{1,2}_{\rm loc}(\Omega)$, called the \emph{stream function} of $u$, such that one has
\begin{equation}\label{e2.3}
\begin{array}{lllll}
 \nabla \tilde u= J \sigma \nabla u &\hbox{in}&\Omega&\ ,&
J:=\left(
\begin{array}{ccc}
0&-1\\
1&0
\end{array}
\right)\ .
\end{array}
\end{equation}
Setting
\begin{equation}
\label{e2.4}
F=u+i \tilde u
\end{equation}
one has
$F=u+i \tilde u \in W^{1,2}_{\rm loc}(\Omega;\mathbb R^2)$ and one writes,  in complex notations,
\begin{equation}\label{1storder}
\begin{array}{ll}
F_{\bar{z}}=\mu F_z +\nu \bar{F_z}\ , & \hbox{in $\Omega$}\ ,
\end{array}
\end{equation}
where, the so called complex dilatations $\mu , \nu$ are given by
\begin{equation}\label{SNU}
\begin{array}{llll}
\mu=\frac{\sigma_{22}-\sigma_{11}-i(\sigma_{12}+\sigma_{21})}{1+{\rm
Tr\,}\sigma +\det \sigma}& \ ,&\nu =\frac{1-\det \sigma
+i(\sigma_{12}-\sigma_{21})}{1+{\rm Tr\,}\sigma +\det \sigma}\ ,
\end{array}
\end{equation}
and satisfy \eqref{ellQC} for some $K\geq 1$ only depending on
$\alpha, \beta$, or in other words $F$ is a quasiregular mapping.

In this paper we are interested in the opposite route, as well.
Given measurable complex valued functions $\mu$ and $\nu$
satisfying \eqref{ellQC}, consider the matrix $\sigma$ defined as
follows
\begin{equation}\label{S}
\sigma:=
\left(
\begin{array}{lll}
\frac{|1-\mu|^2-|\nu|^2}{|1+\nu|^2-|\mu|^2}&\frac{2 \mathfrak{Im}  (\nu -\mu)}{|1+\nu|^2-|\mu|^2} \\
\\
\frac{-2 \mathfrak{Im}  (\nu
+\mu)}{|1+\nu|^2-|\mu|^2}& \frac{|1+\mu|^2-|\nu|^2}{|1+\nu|^2-|\mu|^2}
\end{array}
\right) \ ,
\end{equation}
which is obtained just by inverting the algebraic system
\eqref{SNU}. One can check \cite{aron} that if \eqref{ellQC} holds
for some for given $K\geq 1$, then there exists $\alpha,\beta > 0
$ such that \eqref{ellTM} holds for $\sigma$ as defined in
\eqref{S}. In short, ellipticity in the Beltrami sense implies
ellipticity in the Murat \& Tartar sense.

The exact relationship between $K$ and $(\alpha,\beta)$
will not play a crucial role here. However, we shall prove
the following.

\begin{proposition}\label{ellEG}
Let $(\mu, \nu)$ satisfy the ellipticity condition \eqref{ellQC},
let $\sigma$ be defined via \eqref{S}.  Then $\sigma$ satisfies \eqref{ellTM} with
\begin{equation}\label{alfabeta(K)}
\begin{array}{lll}
\alpha= \frac{1}{K}&\hbox{and}&\beta =K\ .
\end{array}
\end{equation}
Conversely assume that $\sigma\in \mathcal M
(\lambda,\frac{1}{\lambda},\Omega)$ for some $\lambda \in (0,1]$
and let $(\mu, \nu)$ be defined by \eqref{SNU}. Then $(\mu, \nu)$
satisfy the ellipticity condition \eqref{ellQC} with   $K$ defined
as follows
\begin{equation}\label{Koptimal}
K=\frac{1+\sqrt{1- \lambda^2}}{\lambda} \ .
\end{equation}
\end{proposition}
See Section \ref{dimo} for a proof, which also shows the
optimality of these bounds.

\subsection{Quasiconformal solutions to \eqref{e2.3}.}
A question that is crucial in the mere formulation of Conjecture \ref{conj1} is the following.

\noindent
\emph{Is it possible to prescribe a Dirichlet boundary data $g$ on
the real part of $F$ as defined in  \eqref{e2.4} so that the solution to \eqref{1storder} with
that boundary data is globally one-to-one?}

 Or, equivalently, for $\sigma \in {\mathcal M}(\alpha,\beta,\Omega)$, consider the Dirichlet problem
\begin{equation}\label{2nddir}
\left\{
\begin{array}{ll}
{\rm div}(\sigma \nabla u )=0 \ ,& \hbox{in $\Omega$} \ ,\\
u=g\ , & \hbox{on $\partial \Omega$}\ .
\end{array}
\right.
\end{equation}

\noindent
\emph{Under which condition on $g$ the mapping $F= u + i \tilde u$
is one-to-one?}

We recall that solutions to the Beltrami equation \eqref{1storder} are
$K$-quasiregular mapping, therefore the question can be rephrased
as requiring a boundary data which give rise to a global
quasiconformal solution.

Such issues turned out to be very important in  applications of
very different character \cite{ADV,aron,ln,anII,laszlo} and were
addressed already in past years.

The relevant notion in this context is \emph{unimodality}. Assume
that $\partial \Omega$ is a simple closed curve. We say that a
continuous, real valued function $g$ on $\partial \Omega$ is
\emph{unimodal} if $\partial \Omega$ can be split into two simple
arcs on which $g$ is separately monotone (increasing on one arc
and decreasing on the other, once the orientation on $\partial
\Omega$ is fixed). We shall also say that $g$ is \emph{strictly
unimodal} if it is strictly monotone on the same arcs. We shall
prove the following.

\begin{theorem}\label{Thm:onetoone}Let $F\in W^{1,2}_{loc}(\Omega,\mathbb{C})$ be a
solution to \eqref{1storder} such that $u=\mathfrak{Re} F \in
C(\overline{\Omega})$. If $g=u|_{\partial \Omega}$ is unimodal
then $F$ is one-to-one in $\Omega$.
\end{theorem}

The above statement summarizes a circle of reasonings which, in
the last two decades, has been repeatedly used in various contexts
\cite{AannFI,ADV,aron,an}.  See in particular \cite[Proposition
3.7]{aron}, where indeed an interior H\"{o}lder bound for $F^{-1}$
is obtained. A sketch of a proof is given, for the convenience of
the reader in Section \ref{dimo}.

The first result in this direction we are aware of is due to
Leonetti and Nesi \cite[Theorem 5]{ln}. Indeed they proved a
stronger statement.

\noindent\emph{If $g$ is strictly unimodal and $F \in
C(\overline{\Omega};\mathbb{C})$ then $F$ is one-to-one in
$\overline{\Omega}$}.

In fact, in \cite{ln} there are two additional assumptions, that
$\Omega$ is a disk, and that $\sigma$ is symmetric, that is, in
other words, $\mathfrak{Im}\, \nu = 0$. However, such assumptions
are indeed immaterial, in fact we can always reduce to the case
that $\Omega$ is a disk by a conformal mapping, and if $F$ solves
\eqref{1storder} then, as is well-known, it also solves a similar
equation with $\nu=0$ and $\mu$ replaced by
\begin{equation}
\tilde{\mu} = \mu + \frac{\overline{F_z}}{F_z}\nu \ .
\end{equation}

Later, a result of the same sort was proven also in \cite[Theorem
6.1]{B}. In this case the assumptions are that $F\in
W^{1,2}(\Omega,\mathbb{C})$ and that $g=\mathfrak{Re} F_0$ where
$F_0$ is a given quasiconformal mapping whose one-to-one image is
a convex domain. It is worth noticing that this last set of
hypotheses clearly \emph{implies} both $F \in
C(\overline{\Omega};\mathbb{C})$ and the unimodality of $g$.

\subsection{$\sigma$-harmonic mappings.}
Now we review several known results about the so-called
$\sigma$-harmonic mappings. We close this subsection by
reformulating Conjecture \ref{conj1} in the language of
$\sigma$-harmonic mappings and stating Theorem \ref{conj2} which
proves Conjecture \ref{conj1}. Possibly because of a slightly
different language, several results which were published before
\cite{B,B0} may have escaped the authors' attention. We review
here those of more immediate relevance for Conjecture \ref{conj1}
and  postpone a few of them to the following Sections. In order to
rephrase what is already known it is convenient to use the
following notation. We fix $\sigma\in {\mathcal
M}(\alpha,\beta,\Omega)$ and we denote by $U=(u_1,u_2)$ the
$W^{1,2}(\Omega,\mathbb R^2)$ solution to
\begin{equation}\label{pair2ndorder}
 \left\{
 \begin{array}{lc}
{\rm div}(\sigma \nabla u_1)=0 \ ,& \hbox{in $\Omega$}\ ,\\
u_1=x_1\ ,& \hbox{on $\partial \Omega$}\ ,\\
{\rm div}(\sigma \nabla u_2 )=0 \ ,& \hbox{in $\Omega$}\ ,\\
u_2=x_2\ , & \hbox{on $\partial \Omega$} \ .\\
\end{array}
\right.
\end{equation}
Finally we define the stream functions associated to $u_1$ and $u_2$ to be $\tilde u_1$ and $\tilde u_2$ respectively.
Using these notations and recalling (\ref{pair1storder}), we have the identities
\begin{equation}\label{translation}
\begin{array}{lll}
\Phi \equiv u_1 +i \tilde u_1&\ ,& \Psi \equiv u_2 + i \tilde u_2\
.
\end{array}
\end{equation}
Alessandrini and Nesi use the terms   $\sigma$-harmonic functions and  $\sigma$-harmonic mapping
 for $u_1, u_2$ and  $U$ respectively.
With this language, one can compute
\begin{equation}\label{equival}
\mathfrak{Im} (\Phi_z \overline{\Psi_z})= (1+{\rm Tr\,}\sigma
+\det \sigma)\det DU \ .
\end{equation}
Note also that \eqref{ellTM} implies
\begin{equation*}
{\rm Tr\,}\sigma \geq 2\alpha \ , \frac{{\rm Tr\,}\sigma}{\det \sigma} \geq 2\beta^{-1} \ ,
\end{equation*}
and hence
\begin{equation}\label{tracedet}
(1+{\rm Tr\,}\sigma +\det \sigma) > 0 \ .
\end{equation}
The interest of these calculations shall be evident after the following Theorem and Remark.
\begin{theorem}\label{conj2}
Let $\sigma\in {\mathcal M}(K^{-1},K,\Omega)$. If  $\Omega$ is
convex, then the $\sigma$-harmonic mapping $U$ defined by
(\ref{pair2ndorder}) satisfies
 \begin{equation}\label{det>0}
\begin{array}{llll}
\det DU>0&\hbox{almost everywhere}&\hbox{in}&\Omega \ .
\end{array}
\end{equation}
\end{theorem}
\begin{remark}\label{metathm} It is a
straightforward matter to conclude that, by \eqref{equival} and \eqref{tracedet},  Theorem \ref{conj2} proves
Conjecture \ref{conj1} and, consequently, Theorem \ref{main}.
\end{remark}
A proof of Theorem \ref{conj2} will be given in Section
\ref{jacobian}.

The first result towards Theorem \ref{conj2} was proven by Bauman,
Marini and Nesi \cite{bmn}. They proved the assertion under the
assumption that  $\sigma$ is symmetric and of class $C^{\alpha}$.
A further advance was obtained by Alessandrini and Nesi \cite{an}
under the assumption that $\sigma$ is symmetric with  measurable
entries. The two papers follow a common scheme, first one proves
that under suitable conditions on the boundary data (which are
indeed satisfied for the problem \eqref{pair2ndorder} when
$\Omega$ is convex) the mapping $U$ is one-to-one. Here the
guiding light is a conjecture by Rad\`{o} \cite{r}, which was
first proved by H.Kneser \cite{k}  and later, independently, by
Choquet \cite{c}, in the case when $U$ is harmonic. See Theorem
\ref{Kneser} below, for further details. Second, one proves that
if $U$ is locally injective, and sense preserving, then $\det DU>
0$ almost everywhere. In this case the paradigmatic result, in the
harmonic setting, is due to H. Lewy \cite{le}. Actually, in the
harmonic case, and in the case $\sigma\in C^{\alpha}$, one obtains
that $\det DU$ is strictly positive, uniformly on compact subsets.
In the case when $\sigma$ has measurable entries, such uniform
bound cannot hold true. Instead, in \cite{an} it is proven that
for any subset $D$ compactly contained in $\Omega$ one has
\begin{equation}\label{BMO}
\log (\det DU) \in \textrm{BMO}(D)
\end{equation}
which, as is well-known implies that  there  exist  $C\ , \epsilon > 0$ such that in any square $Q\subset \Omega$ one has
\begin{equation}\label{negint}
\left(\frac{1}{|Q|} \int_Q (\det DU)^{\epsilon}dx\right)
\left(\frac{1}{|Q|} \int_Q (\det DU)^{-\epsilon}dx\right) \leq C \
\end{equation}
which clearly implies Theorem \ref{conj2}.

Therefore, when $\sigma$ is symmetric, the tools to prove Conjecture \ref{conj1} were
already available. Later Bojarski, D'Onofrio, Iwaniec and
Sbordone addressed the more general question in the case when
$\sigma$ is not necessarily symmetric. They proved Conjecture
\ref{conj1} in two cases. First when the coefficients are
H\"{o}lder continuous so extending the results by Bauman et al. to
the non-symmetric case. Second they proved the result when $K\leq
3$ so extending the result of Alessandrini and Nesi to the
non-symmetric case in that regime.

In the next two Sections we shall show that the procedure outlined
above for the symmetric case and developed by the authors in
\cite{an} also apply to the non-symmetric case. In fact these
proofs  already appeared  in 2003 as a part of the \emph{Laurea
Thesis} of Natascia Fumolo \cite{nf}, an undergraduate student of
the first author. In this paper we present a much shorter version
by outlining the very few slight changes needed to adapt the
arguments in \cite{an}. On the other hand, some more delicate
issues concerning the precise ellipticity constants, like in
Proposition \ref{ellEG} are treated in a more efficient way here.

In Section \ref{basic} below, we summarize some of the results
obtained in \cite{an} which extend to the non-symmetric case in a
straightforward fashion.

Section \ref{jacobian} contains the core results of this paper,
the main result being Theorem~\ref{t3.1}. From the standpoint of primary pairs
the main implication is Corollary \ref{pr-sol}.

In Section \ref{periodic} we discuss consequences and improvements
to Theorem~\ref{t3.1} in the case of periodic conductivities
$\sigma$, which is relevant in the context of homogenization and
also in connection to  issues concerning the \emph{rigidity} of
gradient fields where quasiconvex hulls are defined either by
using affine or periodic boundary conditions. We refer to
\cite{laszlo}, \cite{acn}, \cite{achn}, \cite{a}  for more
details. The main result here is Theorem~\ref{tA}, which provides
a novel, stronger, quantitative formulation of the non-vanishing
of the Jacobian determinant, in terms of Muckenhoupt weights.

Section \ref{dimo} contains proofs of some auxiliary results.

The final Section \ref{app} collects further developments, remarks
and connections with various relevant areas and applications. In
\S \,\ref{area} we extend some area formulas first discussed in
\cite{anL}. In \S \,\ref{correct} we lay a bridge towards the theory
of \emph{correctors} in homogenization. Finally \S \,\ref{AAstala} develops
an application of the Theorem by Astala \cite{astala},
generalizing results in \cite{ln} and \cite{anL}.

\section{Preliminaries.}\label{basic}
In this Section, $\Omega$ is a  simply connected open subset of
$\mathbb R^2$ and, for applications which will be discussed in
Section \ref{periodic}, we also admit here that $\Omega$ be
unbounded,  possibly  the whole $\mathbb R^2$.
We consider matrix valued  functions $\sigma \in  {\mathcal M}(\alpha,\beta,\Omega)$ as defined in (\ref{ellTM}).

 \begin{notation}\label{not1}
Let $\sigma \in {\mathcal M}(\alpha,\beta,\Omega)$ and let
$U=(u_1,u_2)\in W^{1,2}_{\rm loc}(\Omega,\mathbb R^2)$ be
$\sigma$-harmonic. We denote by $\tilde{U}:=(\tilde u_1,\tilde
u_2)$  the \emph{vectorial stream function} associated to $U$.
Moreover, for any given non zero constant vector $\xi$ we set
$f=U\cdot \xi +i\, \tilde{U}\cdot \xi$.
\end{notation}
\begin{proposition}
Let $\Omega\subseteq \mathbb R^2$ be simply connected and open.
Let $\sigma \in {\mathcal M}(\alpha,\beta,\Omega)$ and let
$U=(u_1,u_2)\in W^{1,2}_{\rm loc}(\Omega,\mathbb R^2)$ be
$\sigma$-harmonic. If for every non zero $\xi$, $f$ is univalent,
then $U$ is univalent.
\end{proposition}
The proof is identical to the proof of Proposition 1 in \cite{an}.
In the latter symmetry of $\sigma$ was assumed but never used.
Details can be found in \cite{nf}.

\begin{theorem}\label{t2.2}
Let $\Omega\subseteq \mathbb R^2$ be a  simply connected and open
set. Let $\sigma \in {\mathcal M}(\alpha,\beta,\Omega)$ and let
$U=(u_1,u_2)\in W^{1,2}_{\rm loc}(\Omega,\mathbb R^2)$ be
$\sigma$-harmonic. Adopt the Notation \ref{not1}. We have that the
following properties are equivalent:
\begin{equation}
\begin{array}{lll}
(i)&\hbox{$f$ is locally one-to-one for every non zero vector $\xi$}\ ,\\
(ii)&\hbox{$U$ is locally one-to-one for every non zero vector $\xi$}\ ,\\
(iii)&\hbox{$\tilde{U}$ is locally one-to-one for every non zero
vector $\xi$}\ .
\end{array}
\end{equation}
\end{theorem}
Also in this case, the proof is identical to the proof of Theorem
3 in \cite{an}, since symmetry of $\sigma$ was assumed but never
used. In fact, additional equivalent conditions to $(i)-(iii)$
were stated in \cite{an}, which involve the notion of
\emph{geometrical critical point}, we omit them here for the sake
of simplicity. Details can be found in \cite{nf}.

\begin{theorem}\label{Kneser}
Let $\Omega$ be a  bounded open set whose boundary  is a simple
closed curve and let $\sigma\in {\mathcal
M}(\alpha,\beta,\Omega)$. Let $\phi=(\phi_1,\phi_2):\partial
\Omega\to \mathbb R^2$ be a sense preserving homeomorphism of
$\partial \Omega$ onto a simple closed  curve $\Gamma$ which is
the boundary of a \emph{convex} domain  $D$. Let $U\in
W^{1,2}_{\rm loc}(\Omega;\mathbb R^2)\cap
C^0(\overline{\Omega};\mathbb R^2)$ be the $\sigma$-harmonic
mapping with components $u_1$ and $u_2$ solving
\begin{equation}
\left\{
\begin{array}{lccc}
{\rm div}(\sigma(x) \nabla u_i(x) )=0 \ ,& \hbox{in $\Omega$}&i=1,2\ ,\\
u_i=\phi_i\ ,& \hbox{on $\partial \Omega$}&i=1,2\ .\\
\end{array}
\right.
\end{equation}
Then
\begin{equation}
\hbox{ $U$ is a sense preserving homeomorphism of
$\overline{\Omega}$ onto $\overline{D}$}\ .
\end{equation}
\end{theorem}
Again, the proof is identical to the proof of Theorem 4 in
\cite{an}, and details can be found in \cite{nf}. Theorem
\ref{Kneser}  generalizes to the measurable, non-symmetric,
context the celebrated result of H. Kneser \cite{k} who solved a
problem raised by Rad\`{o} \cite{r}.

\section{Jacobian of a $\sigma$-harmonic mapping: the \textrm{BMO} bound.}\label{jacobian}
The main subject of this Section is the proof of  Theorem
\ref{conj2}. We will preliminarily proof a much more general
result, namely Theorem \ref{t3.1}.

We recall that, given an  open set $D\subset \mathbb R^2$, $\phi \in L^1_{\rm loc} (D)$ belongs to  ${\rm BMO} (D)$ if
\begin{equation*}
\| \phi \| _{*} = \sup_{Q\subset D} \left( \frac{1}{\mid Q\mid}
\int_Q \mid \phi - \phi_Q\mid \right) <\infty
\end{equation*}
where $Q$ is any square in $D$ and $\phi_Q= \frac{1}{\mid Q\mid} \int_Q \phi$. Recall also that the normed space
 $({\rm BMO} (D), \|\cdot\|_*)$ is in fact a Banach space. The main object of this Section is the following.
\begin{theorem}
\label{t3.1}
Let $\Omega$ be an open subset of $\mathbb R^2$, let $\sigma \in {\mathcal M}(\alpha,\beta,\Omega)$  and
let $U \in W^{1,2}_{\rm loc}(\Omega,\mathbb R^2)$ be a $\sigma$-harmonic mapping  which is locally one-to-one and sense preserving. For every $D\subset \subset \Omega$ we have
\begin{equation}
\label{e3.1}
\begin{array}{ll}
\log (\det DU) \in {\rm BMO} (D)&.
\end{array}
\end{equation}
\end{theorem}
\begin{corollary}\label{pr-sol}
Let $(\mu,\nu)$ be a Beltrami
pair satisfying (\ref{ellQC}) and let $\Phi$ and $\Psi$ be the
solutions to
 \begin{equation}
 \left\{
\begin{array}{lr}
\Phi_{\bar{z}}=\mu \Phi_z +\nu \overline{\Phi_z}\ , & \hbox{in $\Omega$}\ ,\\
\mathfrak{Re} \Phi = \phi_1 \ , & \hbox{on $\partial \Omega$}\ , \\
\Psi_{\bar{z}}=\mu \Psi_z +\nu \overline{\Psi_z}\ ,& \hbox{in $\Omega$}\ ,\\
\mathfrak{Re} \Psi = \phi_2\ , & \hbox{on $\partial \Omega$} \ ,
\end{array}
\right.
\end{equation}
where $\phi=(\phi_1,\phi_2)$, as in Theorem \ref{Kneser}, defines the convex set $D$.
Then  $\Phi$ and $\Psi$ are quasiconformal mappings defined on $\Omega$ which satisfy the inequality
\begin{equation}
\begin{array}{llll}
\mathfrak{Im} (\Phi_z \overline{\Psi_z})>0&\hbox{almost everywhere}&\hbox{in}&\Omega.
\end{array}
\end{equation}
\end{corollary}
\medskip
The proof of Theorem~\ref{t3.1} needs some preparation. It will be
presented at the end of this Section. This part requires slightly
more extended changes with respect to the work in \cite{an}. For
this reason more details will be given.

We recall below two fundamental results, Theorems~\ref{t3.2} and \ref{t3.3}, which will be needed for a proof of Theorem~\ref{t3.1}.

\begin{theorem}[Reimann \cite{re}]\label{t3.2}
Let $f$ be a quasiregular mapping on the open set $D\subset \mathbb R^2$, then for every $D^{\prime} \subset \subset D$
\begin{equation*}
\log (\det Df) \in {\rm BMO} (D^{\prime}) \ .
\end{equation*}
\end{theorem}

\noindent
{\bf Proof.}
See \cite[Theorem 1, Remark 2]{re}. $\Box$

\begin{theorem}[Reimann \cite{re}]
\label{t3.3}
Let $f:D\to G$ be a quasiconformal mapping, $D,G\subset \mathbb R^2$. For every $D^{\prime} \subset \subset D$, there exists $C>0$ such that
\begin{equation*}
\begin{array}{llllll}
\|v \circ f \|_* \leq C  \| v \|_*\ , &{\rm for}&{\rm every}&v\in
{\rm BMO}(f(D^{\prime})) &.
\end{array}
\end{equation*}
\end{theorem}

\noindent
{\bf Proof.} See \cite[Theorem 4]{re} and also \cite[p. 58]{jones}.$\Box$

The next Theorem requires the notion of adjoint equation for a
nondivergence elliptic operator. Let $G \subset \mathbb R^2$ be an
open set. Let $a\in {\mathcal M}(\alpha,\beta,G)$. Set
\begin{equation*}
\begin{array}{ll}
L = \sum\limits_{i,j=1}^2 a_{ij}\frac{\partial^2}{\partial x_i
\partial x_j}&\ .
\end{array}
\end{equation*}
We say that $v\in L^1_{\rm loc} (G)$ is a weak solution of the adjoint equation
\begin{equation}
\label{e3.6}
\begin{array}{lll}
L^* v = 0 \ ,&{\rm in} &G \ ,
\end{array}
\end{equation}
if
\begin{equation*}
\begin{array}{lllll}
\int_G v L u =0\ ,&{\rm for}&{\rm every}& u\in W_0^{2,2}(G)&.
\end{array}
\end{equation*}
We remark that, usually, the ellipticity bounds for  $a$ are
expressed in the form (\ref{ellGT}), rather than (\ref{ellTM}),
but this plays no role here.

\begin{theorem}[Bauman \cite{bau} and Fabes \& Strook \cite{fs}]
\label{t3.4}
For every $w\in L^2_{\rm loc}(G)$, $w\geq 0$, which is a weak solution of the adjoint equation (\ref{e3.6}) we have
\begin{equation}
\label{e3.7}
\left(\frac{1}{\mid Q \mid}  \int_Q w^2\right)^{\frac{1}{2}} \leq  C\left(
\frac{1}{\mid Q \mid}  \int_Q w\right)
\end{equation}
for every square $Q$ such that $ 2 Q \subset G$. Here $C>0$ only depends on the ellipticity constants $\alpha$ and $\beta$.
\end{theorem}

\noindent{\bf Proof.} This Theorem is a slight adaptation between
\cite[Theorem 3.3]{bau} and \cite[Theorem 2.1]{fs}. A proof is
readily obtained by following the arguments in \cite{fs}. The only
additional ingredient which is needed here, is the observation
that, with no need of any smoothness assumption on the
coefficients of $L$, for the special case when the dimension is
two (which is of interest here), for any ball $B \subset G$ and
any $f \in L^2(B)$ there exists and it is unique, the strong
solution
\begin{equation*}
u \in W^{2,2}(B)\cap W_0^{1,2}(B)
\end{equation*}
to the Dirichlet problem
\begin{equation*} \left \{
\begin{array}{lll}
Lu = f \ , & {\rm in}& B\  ,
\\
u=0 \ , &{\rm on} &\partial B\ ,
\end{array}
\right. \end{equation*}
see \cite[Theorem 3]{ta}. $\Box$

\noindent
\textbf{Proof of Theorem~\ref{t3.1}. Preparation.} Let $U= (u_1,u_2)$ satisfy the hypotheses of Theorem~\ref{t3.1} and let
\begin{equation}
\label{e3.8} f = u_1+ i \tilde{u}_1
\end{equation}
be the quasiregular mapping  introduced in Notation~\ref{not1}
with $\xi=(1,0)$. In view of Theorem~\ref{t2.2}, for every $z\in
\Omega$, we can find a neighborhood $D$ of $z$, $D\subset \subset
\Omega$ such that  $U|_D$ and $f|_D$ (i.e. the restrictions of $U$
and $f$ to $D$) are univalent. Therefore, for the proof of
Theorem~\ref{t3.1}, it suffices to show that (\ref{e3.1}) holds
for any sufficiently small $D \subset \subset \Omega$, such that
$U|_D$ and $f|_D$ are univalent. We set
\begin{equation*}
G = f|_D(D)
\end{equation*}
and $V:G\to \mathbb R^2$ given by
\begin{equation}
\label{e3.9}
V= U|_D \circ (f|_D)^{-1}
\end{equation}
where, by definition $(f|_D)^{-1}:G\to D$.  From now on, with a
slight abuse of notation, we will drop the subscripts denoting
restrictions to $D$. We have $DU = (DV \circ f) Df$, and hence
\begin{equation}
\label{e3.10}
\log (\det DU) = \log (\det DV) \circ f + \log (\det Df)~~.
\end{equation}
In view of Theorems~\ref{t3.2} and \ref{t3.3}, the thesis will be
proven as soon as we show that $\log (\det DV)$ belongs to
\textrm{BMO} on compact subsets of $G$. The advantage in replacing
$U$ by $V$, lies in the observation that, in contrast with $\det
DU$,  $\det DV$ satisfies an equation of the type (\ref{e3.6}) for
a suitable choice of the operator $L^*$.

In fact, letting $v_1$ and $\tilde{v}_1$ be the first component of $V$ and its stream function respectively, we can compute
\begin{equation*}
v_1(z) = u_1 \circ f^{-1}(z) = u_1 \circ (u_1+ i
\tilde{u}_1)^{-1}(z) = x_1\ ,
\end{equation*}
\begin{equation}
\label{e3.11} \tilde{v}_1(z) = \tilde{u}_1 \circ f^{-1}(z) =
\tilde{u}_1 \circ (u_1+ i \tilde{u}_1)^{-1}(z) = x_2\ .
\end{equation}
Moreover, by definition,
\begin{equation}
\label{e3.12} \nabla \tilde{v}_1 = J \tau \nabla v_1 \ ,
\end{equation}
where
\begin{equation}
\label{e.13} \tau = T_f \sigma=\frac{Df \sigma Df^T}{\det Df}
\circ f^{-1}\ .
\end{equation}
Hence, using (\ref{e3.11}) and (\ref{e3.12})
\begin{equation*}
 \left(
\begin{array}{ll}
0\\
1
\end{array}\right)
= \left(\begin{array}{ll}
0&-1\\
1&0
\end{array}
\right)
 \left(
\begin{array}{ll}
\tau_{11}&\tau_{12}\\
\tau_{12}&\tau_{22}
\end{array}
\right)
\left(
\begin{array}{ll}
 1\\
0
\end{array}
\right)\ ,
\end{equation*}
that is
\begin{equation}
\label{e3.14}
\tau=
 \left(
\begin{array}{ll}
 1&b\\
0&c
\end{array}
\right)
\end{equation}
where, by construction,
\begin{equation}\label{cb}
\begin{array}{lll}
c= &\det \tau &= \det (\sigma \circ f^{-1})\in L^{\infty}(G)\ ,\\
b= &\tau_{12}&=(\sigma_{12}-\sigma_{21})\circ f^{-1}\in
L^{\infty}(G)\ .
\end{array}
\end{equation}
For a given  $\sigma$, let us denote
\begin{equation}
\begin{array}{lll}
\alpha_{\sigma}= \hbox{ $ \mathrm{ess}\inf\limits_{z\in \Omega}
\left\{ \sigma(z)\xi\cdot \xi \right. $ such that $ \left.\
\xi\in\mathbb{R}^2, \ |\xi|=1 \right\} $}\ ,
\\
\frac{1}{\beta_{\sigma}}= \hbox{$\mathrm{ess}\inf\limits_{z\in
\Omega}\left\{(\sigma(z))^{-1}\xi\cdot \xi \right.$ such that $
\left. \ \xi\in\mathbb{R}^2,\  |\xi|=1\right\}$}\ ,
\end{array}
\end{equation}
that is, $\alpha_{\sigma}, \beta_{\sigma}$ are the best
ellipticity constants $\alpha, \beta$ for which $\sigma \in
\mathcal M (\alpha,\beta,\Omega)$ holds.  We restrict our
attention to the case when
$\alpha_{\sigma}={\beta_{\sigma}}^{-1}:=K^{-1}$. A calculation
that we omit shows that, if $\alpha_{\tau} , \beta_{\tau}$  are
defined accordingly for $\tau$ in $G$, we have
\begin{equation}
\begin{array}{lll}
\label{e3.15}
\alpha_{\tau}= {\rm ess}\inf\limits_{z\in G} \left\{\frac{c(z)+1 -\sqrt{(c(z)-1)^2+b(z)^2}}{2}\ \right\}\ ,\\
 \frac{1}{\beta_{\tau}}  = {\rm ess}  \inf\limits_{z\in G} \left\{\frac{c(z)+1 -\sqrt{(c(z)-1)^2+b(z)^2}}{2 c(z)}\ \right\}\ .
\end{array}
\end{equation}
That is $\tau$  is elliptic in the sense of (\ref{ellTM}) and a
calculation  shows that, in fact, one can take
\begin{equation}\label{alfabeta(tau)}
\alpha_{\tau}=\frac{1}{\beta_{\tau}} =1-\sqrt{1-\frac{1}{K^2}} \ .
\end{equation}
See Section \ref{dimo} for a proof.
Furthermore, by  \eqref{e3.9} and (\ref{e3.11}),
\begin{equation}\label{integrability}
\det DV = \frac{\partial v_2}{\partial x_2}  \in L^2(G)\ .
\end{equation}
Consequently, $v_2$ satisfies
\begin{equation*}
 \hbox{$ \frac{\partial }{\partial x_1} \left(\frac{\partial v_2}{\partial x_1}+b
 \frac{\partial v_2}{\partial x_2}\right)+\frac{\partial}{\partial
x_2} \left(c\frac{\partial v_2}{\partial x_2}\right)= 0$ \ weakly
in $G$}\ .
\end{equation*}
Differentiating the equation above with respect to $x_2$, we see that $w=\det DV$ is a distributional solution of
\begin{equation*}
\begin{array}{lllll}
\frac{\partial^2 }{\partial x_1^2} w+
\frac{\partial^2 }{\partial x_1 \partial x_2}(b w)+
\frac{\partial^2 }{\partial x_2^2} (c w)
 =0 \ , &{\rm in}& G \ ,
\end{array}
\end{equation*}
that is, it is a distributional solution to the adjoint equation
\begin{equation}
\label{e3.16}
\begin{array}{lllll}
L^* w =0\ ,&{\rm in}& G
\end{array}
\end{equation}
where
\begin{equation*}
L= \frac{\partial^2}{\partial x_1^2} + b
\frac{\partial^2}{\partial x_1 \partial x_2}+c
\frac{\partial^2}{\partial x_2^2}\ .
\end{equation*}
On use of (\ref{e3.16}) and (\ref{e3.15}) we may now  apply Theorem~\ref{t3.4}.

We summarize  the resulting statement below.
\begin{proposition}
\label{prop3.5}
For every square $Q$ such that $2Q\subset G$, we have
\begin{equation}\label{reverse}
\left(
\frac{1}{\mid Q\mid} \int_Q (\det DV)^2
\right) ^{\frac{1}{2}}
\leq C \left(\frac{1}{\mid Q\mid} \int_Q \det DV\right),
\end{equation}
where $C>0$ only depends on $\alpha$ and $\beta$.
\end{proposition}

\noindent \textbf{Proof of Theorem~\ref{t3.1}. Conclusion.} A well
known characterization of \textrm{BMO}  in terms of the reverse H\"{o}lder
inequality (see, for instance, \cite[Theorem 2.11 and Corollary
2.18]{gcrdf} ), shows that Proposition~\ref{prop3.5} implies $\log (\det
DV) \in \textrm{BMO}(G^{\prime})$ for every $G^{\prime} \subset \subset G$.
Thus, possibly after replacing $D$ with $D^{\prime} = f^{-1}
(G^{\prime})$, we have, by (\ref{e3.10}) and Theorems~\ref{t3.2}
and \ref{t3.3} that $\log (\det DU) \in \textrm{BMO}(D)$.$\Box$

\noindent \textbf{Proof of Theorem~\ref{conj2}.} Apply
Theorem~\ref{Kneser} with $\phi_1=x_1 \ , \phi_2=x_2$ and
$D=\Omega$, which, by assumption, is convex. Then use
Theorem~\ref{t3.1} . $\Box$

\begin{remark}
We recall now that, in view of Remark \ref{metathm}, the proof of Theorem~\ref{conj2} concludes also the proof of Conjecture~\ref{conj1} and of Theorem~\ref{main}. The proof of Corollary \ref{pr-sol} is also immediate.
\end{remark}

\section{The periodic case.}\label{periodic}
In the homogenization theory, operators with periodic coefficients
play an important role. We refer to the wide literature on the
subject, see for instance \cite{blp} and \cite{m}. We want to remark here that our
result has two interesting consequences in that particular
setting. We set $Q=(0,1)\times (0,1)$  and we shall deal with
functions which are $1$-periodic with respect to each of its
variables $x$ and $y$, which we will call $Q$-{\it periodic}, or
for short, {\it periodic}. For a given $2 \times 2$ matrix $A$, we
write  $U\in W^{1,2}_{\sharp,A}(Q;\mathbb R^2)$ for the space of
zero average (on $Q$) vector fields $U$ such that $U-Ax\in
W^{1,2}_{\sharp}(Q;\mathbb R^2)$, where $W^{1,2}_{\sharp}(Q;\mathbb R^2)$
denotes the completion of $Q$-periodic function with respect to the $W^{1,2}$ norm (see \cite{d} for more details).

We are especially interested in  boundary conditions of periodic
type because of their central
 role in homogenization and in
particular in the so-called $G$-closure problems. In fact, our
starting point for this investigation has its origin in such type
of applications. Given a $2 \times 2$ matrix $ A$, we denote by $U^A=(u_1^A,u_2^A)$
a solution (unique because of our normalization) of
\begin{equation}
\label{e1.4}
\left\{
\begin{array}{lll}
 {\rm div} (\sigma  \nabla u_1^A )= 0\ , \hbox{ in $  \mathbb R^2$\ ,} \\
{\rm div} (\sigma  \nabla u_2^A )= 0\ ,   \hbox{ in $ \mathbb R^2$\ ,}    \\
U^A   \in  W^{1,2}_{\sharp,A} (\mathbb R^2,\mathbb R^2)\ .
\end{array}
\right.
\end{equation}
 The auxiliary problem (\ref{e1.4}) is usually called the cell problem.
 Solutions to (\ref{e1.4}) will be called, with a slight abuse of language, ${\it periodic} \ \sigma$-harmonic mappings.

In the sequel, $\alpha, \beta>0$ and  $\sigma \in {\mathcal M}(\alpha,\beta,\mathbb R^2)$ and $Q$-periodic are given.

\begin{theorem}
\label{tA}
Let $A$ be a non singular $2 \times 2$ matrix  and let $U^A$ be a solution to (\ref{e1.4}). Then we have
\begin{equation}\label{hom}
 \hbox{$U^A$ is a homeomorphism of $\mathbb R^2$  onto itself}.
 \end{equation}
Moreover there exists positive constants $C,\delta$
 only depending on $\alpha$ and $\beta$ such that, for every square $P\subset \mathbb R^2$ and any measurable set $E\subset P$ we have
 \begin{equation} \label{quantitative}
\int_E \frac{\det DU^A}{\det A}\geq
C\left(\frac{|E|}{|P|}\right)^\delta \int_P \frac{\det DU^A}{\det
A}\ .
\end{equation}
\end{theorem}
Here, and in the sequel, integration is meant with respect to
two-dimensional Lebesgue measure.
\begin{remark}
It is worth observing that, when $P=Q$, the unit square, and
$E\subset Q$, we obtain
\begin{equation}\label{global}
\frac{|U^A(E)|}{|\det A|}\geq C |E|^\delta \ .
\end{equation}
Which also  trivially implies
\begin{equation} \label{qualiitative}
\hbox{$\frac{\det DU^A}{\det A} >0 $  almost  everywhere  in $\mathbb R^2\ .$}
\end{equation}
In fact, for any  $\sigma$-harmonic homeomorphism $U$ the area formula
\begin{equation}\label{areaf2}
|U(E)|=\int_E |\det DU|
\end{equation}
holds, see \cite[Proposition 4.2]{anL}, for a proof in the
symmetric case, which however applies equally well to the present
context. See also the discussion in the Section \ref{app} below.
\end{remark}
\begin{remark}
It is anticipated that quantitative
Jacobian bounds, like the one obtained in \eqref{global},  are useful to prove new bounds for effective conductivity i.e. for classes of $H$-limits.
See \cite{n} and \cite{acn}. In particular \cite[Theorem 3.4]{n} gives an explicit improved bound in terms of the constants $C$ and $\delta$ appearing in \eqref{global}. Note the relevance of \eqref{global} in \cite[Definition 3.7]{n}  (thanks to the preceding discussion about the role of the boundary conditions in Section 2 of that paper). However, all such developments would require a careful derivation of bounds for $C$ and $\delta$  and  are beyond
the scope of this note.
\end{remark}

Before beginning the proof  Theorem~\ref{tA}, let us recall some
basic facts about Muckenhoupt weights.

\begin{definition}
A non negative measurable function $w=w(z)$ with $z\in \mathbb C$ is an $A_{\infty}$-weight if

\noindent
(i) there exist constants $C,\delta>0$ such that for every square $P$ and every measurable set $E\subset P$ we have
\begin{equation}\label{5.7}
\frac{\int_E w }{\int_P w }\leq C\left(\frac{|E|}{|P|}\right)^\delta.
\end{equation}
\end{definition}
Thus, as is well-known, the  $A_{\infty}$ condition is a property
of absolute continuity, \emph{uniform  at all scales}, of the
measure $w\textrm{d}x$ with respect to Lebesgue measure
$\textrm{d}x$. The following characterizations of $A_{\infty}$ are
also well-known, see for instance \cite[ Lemma 5]{cf}.
\begin{lemma}\label{Ainfequiv}
Condition (i) above is equivalent to (ii) and (iii) below.

\noindent
(ii) There exist constants $N,\theta>0$ such that for every square $P$
\begin{equation}\label{5.8}
\left(\frac{1}{|P|} \int_P w^{1+\theta} \right)^{\frac{1}{1+\theta}}
\leq N\left(\frac{1}{|P|} \int_P w  \right).
\end{equation}

\noindent
(iii)
There exists constants $M,\eta>0$ such that  for every square $P$ and every measurable set $E\subset P$, we have
\begin{equation}\label{5.9}
\frac{\int_E w }{\int_P w }\geq M\left(\frac{|E|}{|P|}\right)^\eta.
\end{equation}
\end{lemma}
We observe that the quantitative relationships among the pairs of
constants $(C,\delta), (N,\theta)$ and $(M,\eta)$ appearing in the
equivalent characterizations of $A_{\infty}$ can be constructively
evaluated, see Vessella \cite{vess}.

We shall also make use of the following observation.
\begin{remark}\label{sigmatilde}
Let $\sigma\in {\mathcal M}(\alpha,\beta,\Omega)$ and let $u$ be
$\sigma$-harmonic in $\Omega$. Then, up to a multiplicative
scaling, we have that $u$ is also $\tilde{\sigma}$-harmonic with
\begin{equation}\label{eq:sigmatilde}
\tilde{\sigma} = \sqrt{\frac{\beta}{\alpha}}\sigma \in {\mathcal
M}\left(\sqrt{\frac{\alpha}{\beta}},\sqrt{\frac{\beta}{\alpha}},\Omega\right)
\ .
\end{equation}
Thus in the proof below, we may assume, without loss of
generality, $\sigma\in {\mathcal M}(K^{-1},K,\Omega)$ with
$K=\sqrt{\beta / \alpha}$.
\end{remark}

\noindent \textbf{Proof of Theorem~\ref{tA}.}
It suffices to treat
the case when $A$ is the identity matrix $I$ because  $U^A= A
U^{I}$. From now on, for simplicity, we omit the superscript
${I}$.
\smallskip
\noindent The proof of (\ref{hom}) follows with no substantial
changes the one in \cite[Theorem 1]{an}. The proof of
(\ref{quantitative}) consist of showing that $\det DU$ is a
Muckenhoupt  weight. We observe that the arguments of Theorem
\ref{t3.1} tell us that $(\det DU)^\epsilon$ is a Muckenhoupt
weight  for some sufficiently small $\epsilon>0$. Here we improve
the result and show that this is true also for $\epsilon=1$.

By Remark \ref{sigmatilde},  we may assume  $\sigma\in {\mathcal
M}(K^{-1},K,\Omega)$ with $K=\sqrt{\beta / \alpha}$.

Using the notation of Section \ref{jacobian}, we have $U=V \circ
f$ where $f$ now is a $K$-quasiconformal homeomorphism of $\mathbb
C$ onto itself.  Moreover $V$ satisfies (\ref{reverse}) for all
squares in $\mathbb C$. Recall also that $V$ is a $\tau$-harmonic
homeomorphism of  $\mathbb C$ onto itself with $\tau$ given by
(\ref{e3.14}), hence we also have that area formulas of the type
(\ref{areaf2}) also apply to $V$, and obviously to $f$ because of
its quasiconformality.

By (\ref{reverse}) we deduce that $\det DV$ is an $A_{\infty}$-weight, and for suitable $M,\eta>0$ only depending on $K$, we have
\begin{equation}\label{5.10}
\int_F \det DV  \geq M \left(\frac{|F|}{|P|}\right)^{\eta} \int_P
\det DV
\end{equation}
for any square $P$ and any measurable set $F\subset P$.

 Since $f$ is $K$-quasiconformal, we have that
$f$  satisfies the following condition, which can be viewed as one
of the many manifestations of the bounded distortion property of
quasiconformal mappings.

\noindent \emph{There exist $q\in (0,1)$ depending on $K$ only
such that for every square $P\subset \mathbb C$, there exists a
square $P^\prime \subset \mathbb C$ such that}
\begin{equation}
q P^\prime \subset f(P) \subset P^\prime\ .
\end{equation}
Here, if $l$ is the length of the side of $P^\prime$, we denote
by $q P^\prime$ the square concentric to $P^\prime$ with side $q
\cdot l$. We refer to \cite[Proof of Theorem 9.1]{lv} for a proof.

Therefore, we have $f(E)\subset f(P)\subset P^\prime$ and hence
\begin{equation}\label{5.12}
|U(E)| =|V(f(E))|\geq M \left( \frac{|f(E)|}{|P^\prime|}
\right)^\eta |V(P^\prime)|\ .
\end{equation}
Obviously,
\begin{equation*}
\hbox{$|V(P^\prime)|\geq |V(f(P))|$ and $|P^\prime| =\frac{1}{q^2}
|q\,P^\prime| \leq \frac{1}{q^2}|f(P)|$}\ .
\end{equation*}
 Therefore
 \begin{equation}\label{5.13}
 |U(E)| \geq Q \,q^{2\eta} \left(\frac{|f(E)|}{|f(P)|}\right)^\eta |U(P)|\ .
 \end{equation}
By Gehring's Theorem \cite{geh}, we have that $\det Df$ satisfies
a reverse H\"older inequality of the form $(ii)$ in Lemma
\ref{Ainfequiv}, with constants only depending on $K$.
 By $(iii)$ in Lemma
\ref{Ainfequiv}, there exists $L,\rho>0$ only depending on $K$
such that
 \begin{equation}\label{5.14}
 \frac{|f(E)|}{|f(P)|} \geq L \left(\frac{|E|}{|P|}\right)^\rho
 \end{equation}
 and finally, by \eqref{5.13} and \eqref{5.14}
 \begin{equation}
 |U(E)| \geq Q ( q^2\, L)^\eta \left(\frac{|E|}{|P|}\right)^{\eta\,\rho} |U(B)|\ .
 \end{equation}
Thus \eqref{quantitative}  follows.$\Box$

\begin{remark}
The $A_{\infty}$-property of the Jacobian determinant, obtained in
 Theorem~\ref{tA} for the periodic setting,  is indeed an improvement
of the \textrm{BMO} bound obtained previously and which applies to
the wider context of locally injective $\sigma$-harmonic mappings.
Local versions of a bound like \eqref{quantitative} could be
obtained as well for locally injective $\sigma$-harmonic mappings,
however it is expected that a quantitative evaluation of the
constants might be more involved in this case.
\end{remark}

\section{Miscellaneous proofs.}\label{dimo}

\noindent \textbf{Proof of Theorem \ref{Thm:onetoone} (Sketch).} By
the well-known Sto\"{\i}low representation, see for instance \cite[Chapter VI]{lv},
there exists a quasiconformal mapping $\chi: \mathbb{C}\rightarrow
\mathbb{C}$ such that $F$ factorizes as $F=H\circ\chi$ with $H$
holomorphic in $\chi(\Omega)$. Thus, up to the change of variable
$\chi$, one can assume w.l.o.g. $\mu=\nu=0$. Then $u$ is harmonic
and $\tilde{u}$ is its harmonic conjugate. Being $g$ unimodal, $u$
has no critical point inside $\Omega$ \cite{AannFI, amsiam},
moreover, by the maximum principle, for every $t\in(\min g , \max
g)$ the level set $\{u>t\}$ is connected and the level line
$\{u=t\}$ in $\Omega$ is a simple open arc. On $\{u=t\}$,
$\tilde{u}$ has nonzero tangential derivative, hence it is
strictly monotone there. Consequently, $F$ is one-to-one on
$\Omega$.$\Box$

\noindent {\bf Proof of Proposition \ref{ellEG}.} The proof of
this Proposition is a calculus matter regarding matrices $\sigma$
and complex numbers $\mu, \nu$ linked by the relations
\eqref{SNU}, or equivalently \eqref{S}. The dependence on the
space variables $z=x_1+ix_2$ plays no role at this point, and thus
we can neglect it. The inequalities \eqref{ellTM} can be viewed as
lower bounds on the eigenvalues of the symmetric matrices
$\frac{\sigma + \sigma^T}{2}$ and $\frac{\sigma^{-1} +
(\sigma^{-1})^T}{2}$. In terms of $\mu , \nu$, the lower
eigenvalues of such matrices are given by
\begin{equation}
\frac{(1-|\mu|)^2- |\nu|^2}{|1+\nu|^2-|\mu|^2} \ , \ \frac{(1-|\mu|)^2- |\nu|^2}{|1-\nu|^2-|\mu|^2} \ ,
\end{equation}
respectively.
By computing the minima of such expressions as $\mu,\nu \in \mathbb{C}$ satisfy \eqref{ellQC} we obtain \eqref{alfabeta(K)}. It is worth noticing that such minima are achieved when $\nu=|\nu|$ in the first case, and when $\nu=-|\nu|$ in the second case. In either case, the corresponding $\sigma$ turns out to be symmetric.

Viceversa, if we constrain $\mu,\nu$ to satisfy both limitations
\begin{equation}\label{constralfabeta}
\frac{(1-|\mu|)^2- |\nu|^2}{|1+\nu|^2-|\mu|^2}\geq\lambda \ , \ \frac{(1-|\mu|)^2- |\nu|^2}{|1-\nu|^2-|\mu|^2} \geq\lambda\ ,
\end{equation}
then the maximum of $|\mu|+|\nu|$ turns out to be $\sqrt{\frac{1-\lambda}{1+\lambda}}$ and \eqref{Koptimal} follows. Note that in this case the maximum is achieved with $\mu , \nu$ satisfying $ \mu=0$ and $\mathfrak{Re}\nu =0$  which means
\begin{equation}\label{extr-matrices}
\sigma=\left(
\begin{array}{ccc}
a&b\\
-b&a
\end{array}
\right)\hbox{ with } \ a=\lambda \ , \ b = \pm \sqrt{1-\lambda^2}
\ .
\end{equation}
Let us also recall the well-known fact that, if we a-priori assume
$\sigma$ symmetric, then, under the constraints
\eqref{constralfabeta}, the maximum of $|\mu|+|\nu|$  becomes
$\frac{1-\lambda}{1+\lambda}$, that is $K=\frac{1}{\lambda}$.
$\Box$

\noindent\textbf{Proof of \eqref{alfabeta(tau)}.}
As in the Proof of Proposition \ref{ellEG}, we can neglect the
dependence on the space variables $z=x_1+ix_2$. The task here is
to evaluate the minimum eigenvalue of the symmetric part of the
matrices $\tau$ and of $\tau^{-1}$. It suffices to consider the
case $\det \sigma\leq1$. Indeed, up to replacing $\sigma$ with
$\sigma^{-1}$ we can always reduce to this case. Set $D=\det
\sigma$, $T={\rm\,Tr} \sigma$ and $H=(\sigma_{12}-\sigma_{21})^2$.
Elementary computations lead us to minimize the functions
\begin{equation}\label{6a}
F(D,H)=\frac{D+1-\sqrt{(D-1)^2+H}}{2} \ ,
\end{equation}
\begin{equation}\label{6b}
G(F,H)=\frac{F(D,H)}{D}\ ,
\end{equation}
subject to the constraints
\begin{equation}\label{6c}
\frac{T-\sqrt{T^2+H-4D}}{2}\geq \frac{1}{K}\ ,
\end{equation}
\begin{equation}\label{6d}
\frac{T-\sqrt{T^2+H-4D}}{2D}\geq \frac{1}{K}\ .
\end{equation}
Note that, being $D\leq 1$, we have that (\ref{6c}) is always
satisfied if (\ref{6d}) holds and also that $G(D,H)\geq F(D,H)$
with equality when $D=1$. Thus we are reduced to compute
\begin{equation*}
\min\{F(D,H) \, \left|\right. 0\leq D \leq 1 \ , H,T\geq 0 \
,\hbox{(\ref{6c}) holds}\}=1+\sqrt{1-\frac{1}{K^2}}\ .
\end{equation*}
The minimum is achieved when
\begin{equation}\label{minpoint}
\begin{array}{cccc}
T=\frac{2}{K}, &D=1,&\hbox{ and} &H=1-\frac{1}{K^2}
\end{array}
\end{equation}
 which implies
that $\sigma$ has the form (\ref{extr-matrices}) with $\lambda =
1/K$. This proves that $\alpha_\tau$ as defined in \eqref{e3.15}
satisfies \eqref{alfabeta(tau)}. Consequently, by \eqref{6b} and
\eqref{minpoint} we also obtain $
\beta_\tau=\frac{1}{\alpha_\tau}$, proving \eqref{alfabeta(tau)}.
 $\Box$

\section{Further results and connections.}\label{app}
\subsection{Area formulas for $\sigma$-harmonic
mappings.}\label{area}
 One of the original motivations to the
study of Theorem \ref{Thm:onetoone} came from homogenization and
in particular the study of bounds for \emph{effective
conductivity}, that is, $H$-limits. So let $\sigma\in{\mathcal M}
(\alpha,\beta,\mathbb R^2)$ be $Q$-periodic ($Q= (0,1)\times(0,1)
$). By its associated $H$-limit we mean the constant matrix
$\sigma_{\rm eff}$ also called the effective conductivity defined
as the $H$-limit of
$\sigma^{\epsilon}(z):=\sigma(\frac{z}{\epsilon})$ which, as  is
well-known, it is defined via   cell problems as follows. For any
vector $\xi\in \mathbb R^2$, one has
\begin{equation}\label{vp}
\sigma_{\rm eff}\xi \cdot \xi=\min \left\{  \int_Q \sigma \nabla u
\cdot \nabla u  \,\left| \right. u-\xi\cdot x\in
W^{1,2}_{\sharp}(Q;\mathbb R)\right\}.
\end{equation}
Let   $u^\xi$ be the minimizer of (\ref{vp}) and let $\tilde
u^\xi$ be its stream function. Using the notation of Section
\ref{periodic}, we have $u^\xi = U^I \cdot \xi$. Set $f^\xi=u^\xi
+ i \tilde{ u}^{\xi}$. Notice that this quasiconformal mapping
coincides with the one introduced in Notation \ref{not1} when
$U=U^I$. Here we use the superscript $\xi$ just in order to
emphasize this dependence.
\begin{theorem}
For any nonzero vector $\xi\in \mathbb R^2$ one has
\begin{equation}
\sigma_{\rm eff}\xi \cdot \xi=|f^\xi(Q)|.
\end{equation}
\end{theorem}
\noindent {\bf Proof.} We refer to \cite[Proposition 4.1]{anL}.
Again in that context $\sigma$ was assumed to be symmetric but the
hypotheses was not used. $\Box$

The previous result transforms the problem of the calculation of
the effective conductivity into a geometrical one, finding the
area of the set $f^\xi(Q)$.

Next result has already been invoked in  Section \ref{jacobian}.
\begin{theorem}\label{area2}
Let $\Omega$ be a bounded, open, simply connected set. Let
$\sigma\in {\mathcal M}(\alpha,\beta,\Omega)$ and let $U\in
W^{1,2}(\Omega;\mathbb R^2)$ be a  univalent $\sigma$-harmonic
mapping onto an open set $D$. For any measurable set $E\subset
\Omega$ and any function $\phi\in L^1(D;\mathbb R)$ one has
\begin{equation}
\int_E \phi(U(x)) |\det DU(x)| \textrm{d}x  = \int_{U(E)} \phi(y)
\textrm{d}y\ .
\end{equation}
\end{theorem}
\noindent {\bf Proof.} We refer to \cite[Proposition 4.2]{anL}.
Again in that context $\sigma$ was assumed to be symmetric but the
hypotheses was not used. $\Box$

\subsection{Correctors and $H$-convergence.}\label{correct}
In order to explain the meaning of our results in the context of
$H$-convergence we need to recall the notion of correctors. It is
convenient to use the operator  ${\rm Div}$ which acts as the
usual  ${\rm div}$ operator on the rows of $2\times 2$ matrices.

\begin{definition}
Let $\sigma_{\epsilon}$ be a sequence in ${\mathcal M}
(\alpha,\beta,\Omega)$ which is $H$-converging to $\sigma_0$.  Set
$P^{\epsilon}=DU^{\epsilon}$ where, for $\omega$ open with
$\omega\subset \subset \Omega$, one has that $U^{\epsilon}$
satisfies the following properties
\begin{equation}
\left\{
\begin{array}{lrcc}
U^{\epsilon}\in W^{1,2}(\omega;\mathbb R^2)\ ,\\
U^{\epsilon}  \rightharpoonup {\rm Id}\ , &&& \hbox{ weakly in  $W^{1,2}(\omega;\mathbb R^{2 \times 2}) $}\ ,\\
-{\rm Div} ( DU^{\epsilon} (\sigma_{\epsilon})^T) &\to &-{\rm Div}
(\sigma_{0}^T)\ ,   &\hbox{ strongly in  $ W^{-1,2}(\omega;\mathbb
R^2)$\ .}
\end{array}
\right.
\end{equation}
Then $P^\epsilon$ is called a corrector associated with
$(\sigma_{\epsilon},\sigma_0)$.
\end{definition}
For the main properties of the correctors we refer to \cite{mt}.
Let us just recall here that they exist and that, for a given
sequence, $\sigma_{\epsilon}$ which is $H$-converging to
$\sigma_0$,  the difference between two such correctors converges
strongly to zero in $L^2_{\rm loc}(\Omega;\mathbb R^{2\times 2})$.
Our interest in this context is given by the following result.
\begin{proposition}[Murat \& Tartar \cite{mt}]\label{cor-prop}
Let $\sigma_{\epsilon}$  be a sequence in $\mathcal{M}
(\alpha,\beta,\Omega)$ which is $H$-converging to $\sigma_{0}$.
Set $U^{\epsilon}=(u^{\epsilon}_1,u^{\epsilon}_2)\in
H^1(\Omega;\mathbb R^2)$ to be the unique solution to
\begin{equation}\label{correctors}
\left\{
\begin{array}{lll}
{\rm Div} (D U^{\epsilon} \sigma_{\epsilon}^T) ={\rm Div} (\sigma_{0} ^T)\ , &\hbox{ in $\Omega$}\ ,\\
(u_1^{\epsilon},u_2^{\epsilon})= (x_1,x_2)\ , &\hbox{ on $\partial
\Omega$\ .}
\end{array}
\right.
\end{equation}
Then $P^{\epsilon}=DU^{\epsilon}$ is a corrector associated with
$(\sigma_{\epsilon},\sigma_0)$.
\end{proposition}
Proposition \ref{cor-prop} has
a particularly simple interpretation in our language when
$\sigma_{0}$ does not depend on position. In this case (which is
of fundamental importance in the so called $G$-closure problems),
\eqref{correctors} is nothing else than a reformulation of the
boundary value problem \eqref{pair2ndorder}, or equivalently  of
\eqref{pair1storder}, with $\sigma =\sigma_{\epsilon}$ and
Proposition \ref{cor-prop} says that the corrector can be
identified, up to an $L^2$ strong remainder as the Jacobian matrix
of an appropriate $\sigma$-harmonic mapping.

\subsection{Exponent of higher integrability.}\label{AAstala}
As a concluding remark, we observe a straightforward corollary to
Proposition \ref{ellEG} which we state as a Theorem for the
reader's convenience.
\begin{theorem}[Astala]\label{astala}
Let  $\sigma \in {\mathcal M} (\alpha,\beta,\Omega)$  and let
$u\in W^{1,2}_{\rm loc}(\Omega)$ be a $\sigma$-harmonic function.
Set
\begin{equation}\label{kappaastala}
K= \sqrt{\frac{\beta}{\alpha}}+
\sqrt{\frac{\beta-\alpha}{\alpha}}\ .
\end{equation}
Then $u\in W^{1,p}_{\rm loc}(\Omega)$ for any
\begin{equation*}p\in
\left.\left[2,\frac{2K}{K-1}\right. \right)\ .\end{equation*}
\end{theorem}
\noindent\textbf{Proof.} As we noted already in Remark
\ref{sigmatilde}, $u$ is also $\tilde{\sigma}$-harmonic with
$\tilde{\sigma}$ given by \eqref{eq:sigmatilde}, which belongs to
${\mathcal M} (\lambda,\lambda^{-1},\Omega)$ and $\lambda =
\sqrt{\alpha/\beta}$. By Proposition \ref{ellEG}, $f=u+i
\tilde{u}$ is $K$-quasiregular with $K$ given by
\eqref{kappaastala}.  Then one applies the celebrated Astala's
Theorem \cite{astala}.$\Box$

Let us emphasize here that the only, possibly new, observation is of
algebraic nature. In the case when $\sigma$ is symmetric the
{\em algebraically} optimal bound is known as was pointed out in \cite{ln} and
\cite{anL} and achieved for some  $\sigma$'s.
Astala states explicitly in his paper fundamental paper \cite{astala} that the exact
exponent for the $\sigma$-harmonic function seems to depend in a non obvious and
complicated way on the entries of $\sigma$.
Our calculation seems to set
the algebraically optimal  bound in the most general case of non-symmetric $\sigma$. Optimality, in the sense of the existence of a $\sigma$ showing that the exponent of higher integrability cannot be improved, in the context of non
symmetric $\sigma$'s seems to be an open problem. Indeed, by the optimality conditions \eqref{extr-matrices}, the extremal $\sigma$ cannot be symmetric almost everywhere. Therefore it appears that the putative example must be of a new type.

\bigskip
\noindent\textbf{Acknowledgements.} The research of the first
author was supported in part by MiUR, PRIN no. 2006014115 . The
research of the second author was supported in part by MiUR, PRIN
no. 2006017833.

\end{document}